\documentclass[12pt]{article}
\usepackage{enumitem}

\usepackage{graphicx}
\usepackage{amsfonts}
\usepackage{amsmath}
\usepackage{amssymb}
\usepackage[ansinew]{inputenc}
\usepackage{latexsym}
\usepackage{tabularx}
\usepackage{makeidx}
\usepackage{color}
\usepackage{setspace}
\usepackage{multirow}

\allowdisplaybreaks

\begin{document}

\title{\scshape{Hybrid Monte Carlo-Methods in Credit Risk Management}}
\author{Lucia Del Chicca\thanks{supported by IBM Austria}, Gerhard Larcher\thanks{supported by the Austrian Science Fund (FWF): Project F5507-N26, which is a part of the Special Research Program ``Quasi-Monte Carlo Methods: Theory and Applications"}}
\date{}
\maketitle

\onehalfspacing

\begin{abstract}
In this paper we analyze and compare the use of Monte Carlo, Quasi-Monte Carlo and hybrid Monte Carlo-methods in the credit risk management system Credit Metrics by J.P.Morgan. We show that hybrid sequences used for simulations, in a suitable way, in many relevant situations, perform better than pure Monte Carlo and pure Quasi-Monte Carlo methods, and they essentially never perform worse than these methods.   
\end{abstract}


\section{The use of hybrid sequences in simulations}

The use of hybrid sequences in simulation problems in different applications is not new. It was suggested for example by Spanier in \cite{Spanier} for transport problems, by Ökten in \cite{Oekten} for derivative pricing, or by Keller \cite{Keller} in computer graphics. In this article, by a hybrid sequence, we understand a sequence $(x_n)_{n \geq 0} $ in a $d-$dimensional unit-cube $\left[0,1 \right)^{d} $ which is obtained by concatenating $l$ sequences $(y^{(i)}_{n})_{n\geq 0}, \quad i=1, \ldots, l$, each in a $d_i-$dimensional unit cube $\left[0,1\right)^{d_i},$ where $d=d_1+ \ldots + d_l,$ and where these component sequences are different kinds of pseudo-random and/or quasi-random point-sequences. In the last years the investigation of distribution properties of such sequences has become a vivid and challenging field of research.

The most common pseudo-random sequences used to build hybrid sequences are sequences generated by a linear congruence generator or by an inversive congruence generator. The most common quasi-random point-sequences are Halton sequences, Kronecker-sequences, Sobol-sequences or Niederreiter sequences (respectively Hammersley point sets, good-lattice point sets, or digital $(t,m,s)-$nets if it is worked with a pre-determined finite number of points). Just informally said, \textit{quasi-random} point sets are designed such that they show good distribution properties whereas \textit{pseudo-random} point sets show a ``better'' random behavior. We do not go into details here in generation and properties of quasi- and pseudo-random-point sequences or hybrid sequences but refer to the standard literature on these topics like for example \cite{Nieder} or \cite{Pill}, to \cite{Nieder1}, \cite{Nieder2}, \cite{Nieder3}, \cite{Nieder4}, \cite{Winterhof} and to the very informative paper \cite{Rosi} which gives a survey on recent developments concerning hybrid sequences, and to the references given in this paper. 

It is a frequent observation that in $s-$dimensional problems based on a given number $N$ of scenarios (generated by $N$ $s-$dimensional pseudo-random or quasi-random point-vectors) QMC-methods (based on quasi-random point sets) perform significantly better if $s$ is ``rather small'' in dependence on $N$, say if $s<f(N)<N,$ and that MC-methods (based on pseudo-random point sets) perform significantly better if $s$ is ``rather large'' in dependence on $N,$ say if $s> g(N) >> f(N).$ Here the functions $f$ and $g$ essentially depend on the concrete simulation problem. So it seems that in simulation problems where $s>g(N)$ holds, the benefits of QMC-methods cannot be utilized. In some applications, however, it was observed that some, say $s'<s$ of the dimensions of the problem are in some sense of much more importance than the other remaining $s-s'$ dimensions of the problem. If moreover $s' <f(N),$ and $s-s' > g(N)$ then the following approach can give a competitive edge: Instead of a pure MC sequence or a pure QMC sequence for carrying out the simulation, use now a hybrid sequence generated by an $s'-$ dimensional QMC sequence of length $N$ and by a suitable $(s-s')-$ dimensional MC sequence of length $N$. This hybrid sequence is applied in such a way that the ``most important'' $s'$ objects in the simulation problem are treated by the QMC- part, whereas the (-in most cases many-) $s-s'$ objects are treated with the MC part. 

Of course it is necessary, that the hybrid sequence as a whole has good distribution properties (here we can rely on the theoretical research mentioned and cited above). If this is guaranteed then the reasoning is the following: If the influence of the $s'$ dimensions treated by QMC is very strong, then using QMC for these $s'<f(N)$ coordinates improves the simulation results. If the influence of these $s'$ dimensions is not really essential, then the use of QMC for these few $s'$ coordinates does not worsen the simulation (since $s'<f(N)$). Hence, in some cases the use of hybrid sequences in this way should improve the performance of the simulation method, and in any case should not deteriorate the performance of the simulation problem.

In this paper we report on our realization of this use of hybrid sequences to the credit risk management system CreditMetrics by J.P.Morgan. This system will turn out to be in some sense ideal for the above explained approach. 

Section 2 gives a short introduction to the aspects of CreditMetrics necessary to understand the simulation procedure (we give only few details on the background of credit risk management). In Section 3 we describe our simulation procedure and in Section 4 we show the results of our experiments. Section 5 concludes the paper.

\section{The Simulation Part of CreditMetrics}
The system CreditMetrics by J.P. Morgan is completely described in the Technical Document \cite{Morgan}. We will recall here only the aspects which are necessary to understand the simulation goal and the simulation procedure. We will make the same choices as in \cite{Morgan} about the type of credit products used, how to compute their present and future values and the risk horizon considered. For all other details, discussion on these choices and for the background in credit risk management we refer to \cite{Morgan}. 

We consider a portfolio of $s$ credit products. To avoid irrelevant technical details we make the following assumptions and notations: 

\begin{itemize}
	\item each credit product has the form of a bond with annual coupon payments;
	\item the remaining time to maturity of each bond is a multiple of a full year;
	\item the next coupon payment occurs exactly in one year for every bond;
	\item we are interested in computing the risk of the bonds over one year horizon. So we will speak about the ``value of a bond at year-end'' meaning ``at the end of the chosen one year risk horizon''; 
	\item each bond is rated following the Standard\&Poor's terminology (for simplicity we consider only 8 rating classes like in \cite{Morgan}). Each bond ``now'' (at time $0$) belongs to one of the following 7 rating categories: AAA, AA, A, BBB, BB, B, CCC. Due to quality changes of the issuer, the rating class of each bond can change once a year to one of the 8 rating categories: AAA, AA, A, BBB, BB, B, CCC, D ($=$ default);
	\item each bond belongs to a certain seniority class (Senior Secured, Senior Unsecured, Senior Subordinated, Subordinated, Junior Subordinated) which determines the recovery rate in event of a default of the bond;
	\item each credit $K_i$ in the portfolio is completely determined by
	
\begin{itemize}
	\item face value $M_i$ (like in \cite{Morgan} for simplicity we assume only one currency for all bonds in the portfolio),
	\item coupon $C_i$ in percent,
	\item time to maturity $T_{i}$  (we assume $T_{i}$ to be a positive integer, we are immediately after a coupon payment),
	\item rating class $X_i$,
	\item seniority class $R_i;$ 
\end{itemize}

  \item the present value $A^{0}_{i}$ of the credit $K_i$ is given by
  
  $$A^{0}_{i}= \sum^{T_{i}-1}_{l=1} \frac{\frac{C_i}{100} \cdot M_i}{\left(1+\frac{f^{X_i}_{0,l}}{100}\right)^{l}} + \frac{\left(1+ \frac{C_i}{100} \right) \cdot M_{i}}{\left(1+\frac{f^{X_i}_{0,T_i}}{100}\right)^{T_{i}}}$$
  
  where $f^{X_i}_{0,T_{i}}$ are the average forward zero rates in percent for $X_i-$rated bonds with maturity $T_i;$

  \item due to quality changes of the issuer, each bond determined by $(M_i, C_i, T_i, X_i, R_i)$ can change to another rating category, say $Y_i$, at year-end. Then its new value at next year-end is given by
  
\begin{equation}
A^{1}_{i}(Y_i) = \sum^{T_{i}-1}_{l=2} \frac{\frac{C_i}{100} \cdot M_i}{\left(1+\frac{f^{Y_i}_{1,l}}{100}\right)^{l-1}} + \frac{\left(1+\frac{C_i}{100}\right) \cdot M_{i}}{\left(1+\frac{f^{Y_i}_{1,T_i}}{100}\right)^{T_{i}-1}}  \quad \mbox{if} \quad Y_i \neq D
\end{equation}  
  
and by

\begin{equation*}
A^{1}_{i}(D)= \frac{RR_i}{100}\cdot M_i \quad \mbox{if} \quad Y_i = D
\end{equation*}  

where $RR_i$ is a certain given recovery rate which depends on the seniority class $R_i$ and $f^{X_i}_{1,T_i}$ are the one-year forward zero rates in percent calculated from $f^{X_i}_{0,1}$ and $f^{X_i}_{0,T_i}$ by
  
  $$f^{X_i}_{1,T_i}= \left(\frac{(1+f^{X_i}_{0,T_{i}})^{T_{i}}}{(1+f^{X_i}_{0,1})}\right)^{\left(\frac{1}{T_{i}-1}\right)}-1 . $$

Following \cite{Morgan} on page 26, in our simulation experiments the values of $RR_i$ are given by Table \ref{tab:recRates};

\begin{table}[!ht]
\centering
\begin{tabular}{l|c|c}
\hline
\textbf{Seniority Class} & \textbf{Mean (\%)} & \textbf{Standard Deviation (\%)} \\ \hline
Senior Secured & $53.80$ & $26.86$ \\ \hline
Senior Unsecured & $51.13$ & $25.45$ \\ \hline
Senior Subordinated & $38.52$ & $23.81$ \\ \hline
Subordinated & $32.74$ & $20.18$ \\ \hline
Junior Subordinated & $17.09$ & $10.90$ \\ \hline
\end{tabular}
\caption{Recovery rates by seniority class. Source: CreditMetrics-Technical Document}
\label{tab:recRates}
\end{table}

\begin{table}
\centering
\begin{tabular}{l|c|c|c|c|c|c|c|c} \hline
 \multirow{2}{14mm}{\textbf{Initial rating}} &  \multicolumn{8}{c} {\textbf{Rating at year-end (\%)}} \\[0.7ex]   \cline{2-9}
  &  AAA & AA & A & BBB & BB & B & CCC & Default \\ \hline
 AAA & 90.81 & 8.33 & 0.68 & 0.06 & 0.12 & 0 & 0 & 0 \\
 AA & 0.70 & 90.65 & 7.79 & 0.64 & 0.06 & 0.14 & 0.02 & 0 \\ 
 A & 0.09 & 2.27 & 91.05 & 5.52 & 0.74 & 0.26 & 0.01 & 0.06 \\ 
 BBB & 0.02 & 0.33 & 5.95 & 86.93 & 5.30 & 1.17 & 0.12 & 0.18 \\
 BB & 0.03 & 0.14 & 0.67 & 7.73 & 80.53 & 8.84 & 1.00 & 1.06 \\
 B & 0 & 0.11 & 0.24 & 0.43 & 6.48 & 83.46 & 4.07 & 5.20 \\
 CCC & 0.22 & 0 & 0.22 & 1.40 & 2.38 & 11.24 & 64.86 & 19.79 \\ \hline
\end{tabular} 
\caption{Rating - transition probabilities. Source: CreditMetrics-Technical Document}
\label{tab:trans}
\end{table}

\item for each credit $K_i$ its expected value at the year end is given by

\begin{equation}
E(A^1_i)= \sum_{Y_i=AAA, AA, \ldots, CCC, D} p_{X_i\rightarrow Y_i} \cdot A^1_i(Y_i)
\end{equation}

where by $p_{X_i\rightarrow Y_i}$ we denote the probability that a now $X_i-$rated credit is $Y_i-$rated at the year-end. For these rating - transition probabilities we use the values given in Table 2 (see \cite{Morgan}, page 20);
  
  \item following $(2)$ the expected value of the entire credit portfolio at the next year-end is given by
  
 \begin{equation}
 EA^1=\sum^{s}_{i=1}E(A^1_i).
 \end{equation}  
\end{itemize}
 
It is one of the main aims of CreditMetrics to determine the $1^{st}-$percentile level $\theta (A^1)$ of the random variable 

$$A^1= \sum^{s}_{i=1}A^1_{i}(Y_i),$$

\noindent that is the value of the credit portfolio at year-end. Whereas by (3) we have an explicit formula $EA^{1}$ for the expected value of $A^{1}$ at year-end, we do not have an explicit formula for $\theta (A^1),$ and we are reliant on simulation methods. To carry out this simulation, CreditMetrics suggests the following procedure (see \cite{Morgan} again for details):

\begin{itemize}
	\item we assume that it is given a correlation matrix $\mathfrak{A}=(\rho_{ij})_{i,j=1, \ldots, s}$ for the asset values of the credits obligors, i.e., $(\rho_{ij})$ is the correlation between the asset values of credit obligors of credits $K_i$ and $K_j.$ We carry out Cholesky-decomposition of $\mathfrak{A}$ i.e., $\mathfrak{A}=W \cdot W^{T};$
	\item we assume that we are given values $Z^{X}_{Y} \in \mathbb{R}$ for $X=AAA, AA, \ldots, B, CCC$ and $Y=AAA, AA, \ldots, CCC, D$ with 
	
$$- \infty \leq Z^{X}_{D} \leq Z^{X}_{CCC} \leq Z^{X}_{B} \leq \ldots \leq Z^{X}_{A} \leq Z^{X}_{AA} \leq + \infty $$

with the following property:

If the ``normalized asset value return at year end'' $Z$ of the obligor of credit $K_i$ with current rating $X_i$ satisfies

\begin{equation}
\begin{array}{ll}
Z \leq Z^{X_i}_{D} & \mbox{ then } K_i \mbox{ has rating } D \mbox{ at year-end, }\\
Z^{X_i}_{D} < Z \leq Z^{i}_{CCC} &  \mbox{ then } K_i \mbox{ has rating } CCC \mbox{ at year-end,}\\
Z^{X_i}_{CCC} < Z \leq Z^{i}_{B} &  \mbox{ then } K_i \mbox{ has rating } B \mbox{ at year-end,}\\
& \ldots  \\
Z^{X_i}_{AA}  \leq Z &  \mbox{ then } K_i \mbox{ has rating } AAA \mbox{ at year-end.}\\
\end{array}
\end{equation}

(Under certain additional assumptions these values $Z^{X}_{Y}$ can be calculated with the help of the values in Table \ref{tab:trans});

 \item now we generate a sample of $N$ independent $s-$dimensional random-vectors each vector consisting of $s$ independent standard normally distributed random variables. These vectors are multiplied by the matrix $W$ and so we obtain $N$ vectors 
 
$$\xi^i = \left( 
\begin{array}{l}
\xi^{i}_{1}\\
\ldots\\
\xi^{i}_{s}\\
\end{array}
\right), \quad i=1, \ldots, N
$$ 

where $\xi^{i}_{j}$ represents the ``next year normalized asset value'' of the obligor of credit $K_j$ in the $i-$ th random sample. $\xi^{i}_{j}$, through $(4)$ determines the rating of credit $K_j$ at year-end in the $i-$th random sample. After doing so for each credit, finally, with the help of $(1)$ we can determine $A^{1}(\mbox{of sample } i)$, the value of the credit portfolio for the $i-$th scenario. In this way we obtain $N$ possible simulation values for the credit portfolio at year-end.  
 \item If we order these values and denote them by
 
$$\tilde{A}^{1}(1) \leq \tilde{A}^{1}(2) \leq \ldots \leq \tilde{A}^{1}(N), $$

then 

$$\tilde{\theta}_N:= \tilde{A}^{1}\left(\left[\frac{N}{100}\right]\right) $$

gives an estimate for the $1^{st}-$percentile $\theta(A_1)$ of the credit-portfolio value at the year-end;

 \item Just to test the simulation procedure we can compare the value 
 
 $$\tilde{E}A^{1}_{N}:= \frac{1}{N} \sum^{N}_{i=1} A^1(\mbox{sample }i)$$
 
 with the exact value for the expectation $EA^1$ of $A^1$ given in (3);

 \item For the case of portfolios of non-correlated credits we also have an exact value for the variance $VA^1$ of the value of the credit portfolio in one year. Of course it is just the sum of the variances $V(A^1_i)$ of the values $A^1_i$, and 

$$V(A^1_i)= \sum_{Y_i=AAA, AA, \ldots, CCC} p_{X_i\rightarrow Y_i} \cdot (A^1_i(Y_i)-E(A^1_i))^2 + p_{X_i\rightarrow D} \cdot V(RR_i),$$

where $V(RR_i)$ is the square of the standard deviation of the recovery rate, given by Table \ref{tab:recRates}. So in this case we can also simulate the standard deviation $\sigma(A^1)=\sqrt{V(A^1)}$ on the credit portfolio and test the simulation methods by comparing them with the exact value.

\end{itemize}

This is the complete simulation procedure which we will carry out in the next section. In fact we will work in the following with ``normed portfolio values'' rather than with ``absolute portfolio values''. That means that, instead of $A^{1} (\mbox {sample } i), \tilde{A^{1}}(i)$ and $EA^{1}$ we consider

$$\frac{100 \cdot A^{1} (\mbox {sample } i)}{M}, \frac{100 \cdot \tilde{A^{1}}(i)}{M} \mbox{ and } \frac{100 \cdot EA_1}{M},$$ 

\noindent where $M:= \sum^{s}_{i=1} M_i$ is the total face value of the credit portfolio and as $\sigma(A^1)$ we take the standard deviation of the normed values. In the following we will use the above notation $\left(EA^{1}, \tilde{E}A^{1}_{N}, \sigma(A^1), A^{1}(\mbox{sample }i), \tilde{A}^{1}(i), \tilde{\theta}_N, \theta \right)$  for the \textit{normed} versions. 

\section{Simulation experiments with MC, QMC and Hybrid sequences} 
In our simulation experiments we have worked with uncorrelated credit portfolios of different size: the test portfolios will consist of $100, 500,$ or $1.000$ credit products. These credit portfolios are artificially chosen, with parameters that are typical for real-life portfolios. For each size we considered two types of portfolios: ``homogeneous'' portfolios that consist of bonds with ``realistic distributed'' risk profiles and ``inhomogeneous'' portfolios that contain a certain number of ``very high'' risk profiles.
For all portfolios we calculate first the (normed) exact expected value $EA_1$ and its standard deviation $\sigma(A^1)$ of the portfolio at year- end, then we carry out the simulation (according to the above described procedure) to obtain approximated values $\tilde{E}A^1_N$ for $EA^{1}$ and $\tilde{\sigma}(A^1_N)$ for $\sigma(A^1)$. As mentioned in the last section, $\tilde{E}A^{1}_{N}$ and $\tilde{\sigma}(A^{1}_{N})$ are calculated just as a test for the simulation procedure. 
Further we carry out simulation (always according to the above described procedure) to obtain approximated values $\tilde{\theta}_N$ for $\theta.$ Here we have no method to obtain the exact reference value $\theta,$ but we obtain approximate reference values $\tilde{\theta}$ by carrying out Monte Carlo simulation with a very large number of samples.   
For our experiments we used up to $1.000$ sample scenarios for portfolios consisting of $100$ credits, up to $5.000$ scenarios for portfolios consisting of $500$ credits, and up to $10.000$ sample scenarios for portfolios consisting of $1.000$ credits. 
For determining an approximate reference value $\tilde{\theta}$ for $\theta$ we used MC simulation with $50.000$ scenarios in each case.
For each size of portfolios we perform MC simulation, QMC simulation with Niederreiter sequences and simulation with hybrid sequences which were generated by a pseudo-random sequence and by a Niederreiter sequence. Of course in a first step in carrying out the simulation procedure described in Section 2, we have to transform the point sets which are drawn from a uniform distribution in a unit cube to standard normally distributed point sets. We use a standard inversion method. 

For the tests with hybrid sequences in advance we have carried out experiments to determine an adequate quantity $f(N)$ for our type of application. It turned out that for our showcases 

$$
\begin{array}{lll}
\mbox{SC I}: & $s=100, N=1.000$\\
\mbox{SC II}: & $s=500, N=5.000$ \\
\mbox{SC III}: & $s=1.000, N=10.000$\\
\end{array}
$$
QMC with Niederreiter sequences essentially in all experiments gave significantly better results than MC for dimensions $f(1.000)=5, f(5.000)=25, f(10.000)=50.$ So when we carried out  simulations for SC I, SC II, SC III with hybrid sequences, we chose $s' = 5, s' = 25,$ and $s' = 50$ respectively. When using hybrid sequences for our simulations, then in a first step we identify the $s'$ credits of the portfolio with the highest risk profile. As ``risk profile'' of credit $K_i$ we use the following quantity

$$P_i:= \left(1- \frac{RR_i}{100}\right) \cdot M_i \cdot p_{{x_i}\rightarrow D}$$
that is, the average expected loss for credit $K_i$ in case of default, weighted with probability of default. Of course this choice of a risk profile in a general case could be refined by taking also into account the position of $K_i$ in the correlation structure of the portfolio. For the $s'$ credits $K_i$ with the highest values for $P_i$ we use the QMC part of the hybrid sequence. For the remaining $s-s' >> s'$ lower risk credits we use the MC part of the hybrid sequence. Whereas in the ``homogeneous'' test-portfolios the $s'$ most risky credits do not differ significantly from the other credits in the portfolio (concerning their risk profile), we have designed the ``inhomogeneous'' test-portfolios in such a way, that $d \leq s'$ credits in the portfolio show a significant higher risk-profile than the average risk-profile of all credits. 

In the following we give typical results of our simulations, i.e., the graphics of convergence of the different simulation methods for showcases SC I, SC II, SC III for the expected normed portfolio value in one year, for the standard deviation of the portfolio value in one year and for the $1^{st}$ percentile of the normed portfolio value in one year, for a homogeneous and for an inhomogeneous portfolio. So we provide 18 graphics in the following. The exact value (in the case of the percentile the approximate exact value) is always marked with a horizontal line. The results on the MC-simulation are shown in blue color the results of simulations with Niederreiter sequences (QMC) are shown in orange and finally the results of simulation with hybrid sequences are shown in pink.

\section{Numerical results}
As anticipated in Section 3, we will present in this section a typical selection of our results. For each of the three showcases we will show the example of a homogeneous portfolio and of an inhomogeneous portfolio of the same dimension. Each example is described through three figures: the first represents the expected normed portfolio value in one year, the second shows the simulation for the standard deviation of the portfolio value in one year and finally the third shows the $1^{st}$ percentile of the normed portfolio value in one year. So we have a total of six examples and $18$ figures that we list at the end of the article. All the figures are to be read in the same way: the exact values are always marked with a horizontal red line. The MC-simulations are shown with a blue curve, the simulations with Niederreiter sequences (QMC) are shown with an orange curve and the simulation with hybrid sequences are shown with a pink curve.

\begin{figure}[!ht]
  \centering
   \includegraphics[width=0.5\textwidth]{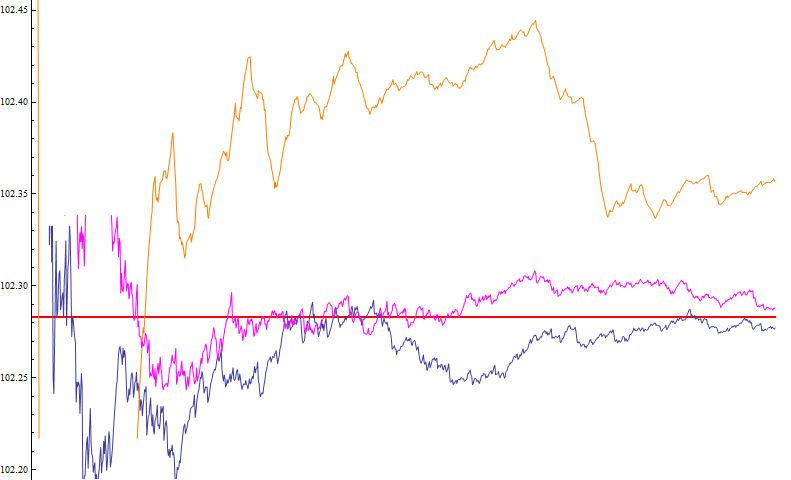}
  \caption{Expected value in the case of a size $100$ homogeneous portfolio and up to $1.000$ simulations, MC (blue line), QMC (orange line) and hybrid simulation (pink).}
  \label{ErwWertHom1000}
\end{figure}

\begin{figure}[!ht]
  \centering
   \includegraphics[width=0.5\textwidth]{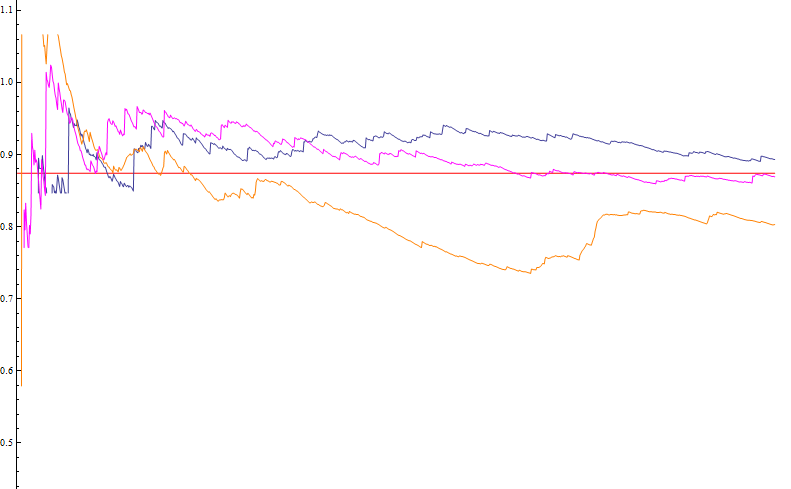}
  \caption{Standard deviation of the portfolio expected value in the case of a size $100$ homogeneous portfolio and up to $1.000$ simulations, MC (blue line), QMC (orange line) and hybrid simulation (pink curve).}
  \label{VolaHom1000}
\end{figure}

\begin{figure}[!ht]
  \centering
   \includegraphics[width=0.5\textwidth]{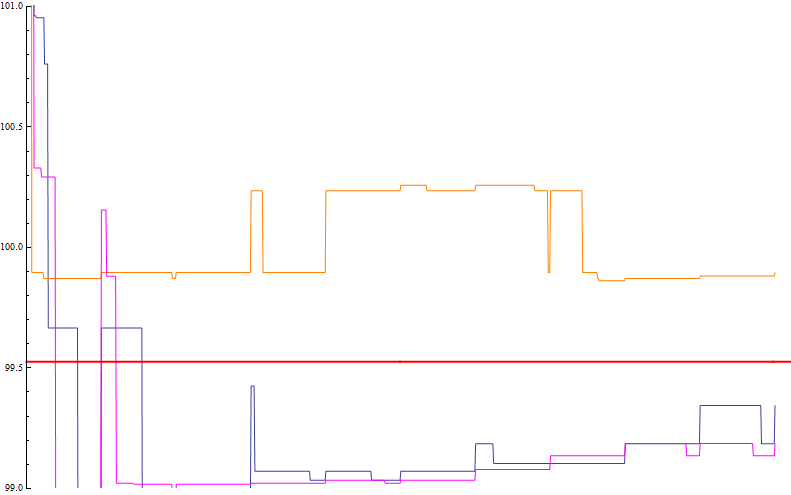}
  \caption{$1^{st}$ percentile level of the portfolio value of a homogeneous portfolio of size $100$ and up to $1.000$ simulations, MC (blue curve), QMC (orange curve) and hybrid simulation (pink curve).}
  \label{PercHom1000}
\end{figure}

The first showcase SC I: $s= 100, N=1.000$ concerns a portfolio of dimension 100 for which we carry out up to 1.000 simulations in the three different methods.
In  Figure  \ref{ErwWertHom1000} we can see that, as expected, the QMC curve oscillates much more than the curves of the other two methods. MC and hybrid method perform very similarly. 
In Figure \ref{VolaHom1000} again MC and hybrid method perform better than QMC method and hybrid slightly better that MC. 
Figure \ref{PercHom1000} shows again a MC and a hybrid method simulation that perform very similarly, both better than QMC method.

\begin{figure}[!ht]
  \centering
   \includegraphics[width=0.5\textwidth]{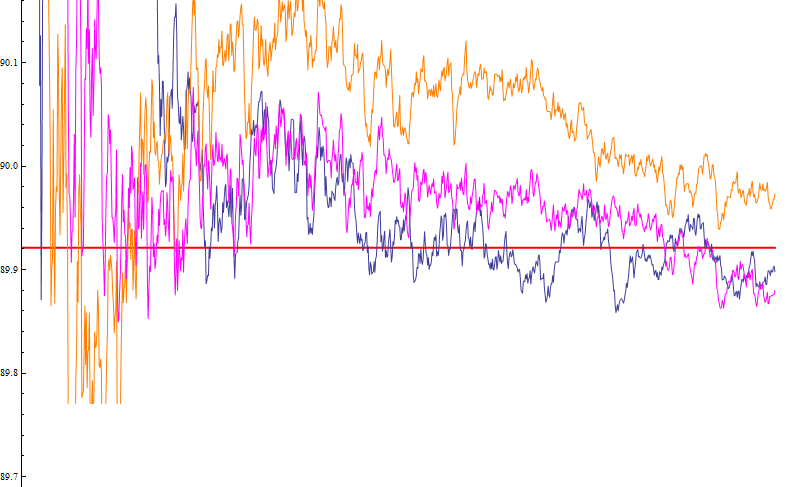}
  \caption{Expected value in the case of a size $100$ inhomogeneous portfolio and up to $1.000$ simulations, MC (blue curve), QMC (orange curve) and hybrid simulation (pink curve).}
  \label{ErwWertInHom1000}
\end{figure}

\begin{figure}[!ht]
  \centering
   \includegraphics[width=0.5\textwidth]{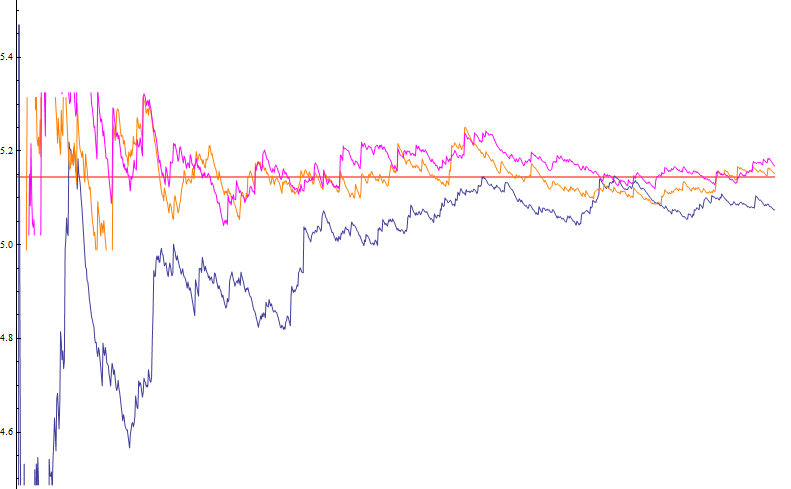}
  \caption{Standard deviation of the portfolio expected value in the case of a size $100$ inhomogeneous portfolio and up to $1.000$ simulations, MC (blue curve), QMC (orange curve) and hybrid simulation (pink curve).}
  \label{VolaInHom1000}
\end{figure}

\begin{figure}[!ht]
  \centering
   \includegraphics[width=0.5\textwidth]{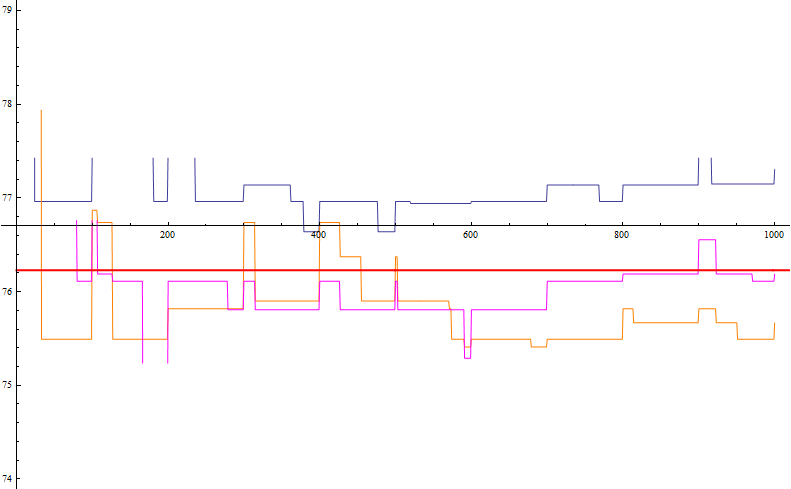}
  \caption{$1^{st}$ percentile level of the portfolio value of a inhomogeneous portfolio of size $100$ and up to $1.000$ simulations, MC (blue curve), QMC (orange curve) and hybrid simulation (pink curve).}
  \label{PercInHom1000}
\end{figure}

We compare these three figures with figures \ref{ErwWertInHom1000}, \ref{VolaInHom1000}, \ref{PercInHom1000} of the example of an inhomogeneous portfolio of the same size. 
In Figure \ref{ErwWertInHom1000} again MC and hybrid methods perform better than QMC. It oscillates much more and even if it converges to the same value it takes longer than the other two methods. 
Figure \ref{VolaInHom1000} shows a much better performance of the hybrid method in comparison with the MC method.
In Figure \ref{PercInHom1000}  it appears immediately that the hybrid method outperforms the other two methods. 

\begin{figure}[!ht]
  \centering
   \includegraphics[width=0.5\textwidth]{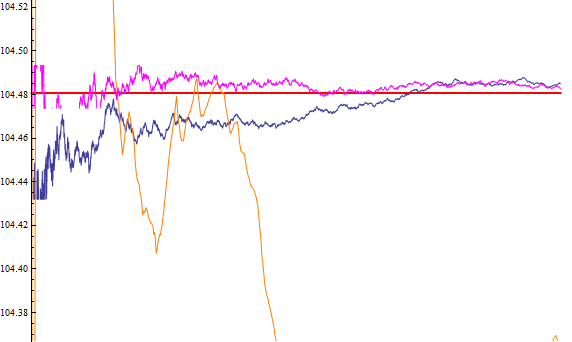}
  \caption{Expected value in the case of a size $500$ homogeneous portfolio and up to $5.000$ simulations, MC (blue curve), QMC (orange curve) and hybrid simulation (pink curve).}
  \label{ErwWertHom5000}
\end{figure}

\begin{figure}[!ht]
  \centering
   \includegraphics[width=0.5\textwidth]{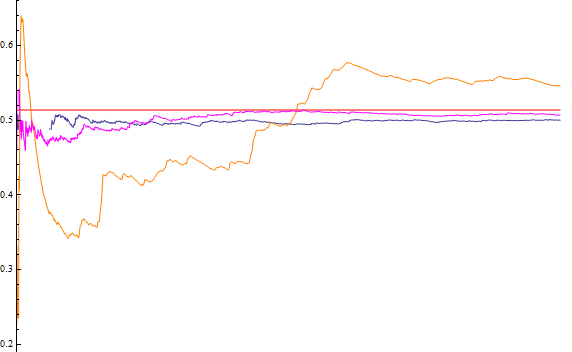}
  \caption{Standard deviation of the portfolio expected value in the case of a size $500$ homogeneous portfolio and up to $5.000$ simulations, MC (blue curve), QMC (orange curve) and hybrid simulation (pink curve).}
  \label{VolaHom5000}
\end{figure}

\begin{figure}[!ht]
  \centering
   \includegraphics[width=0.5\textwidth]{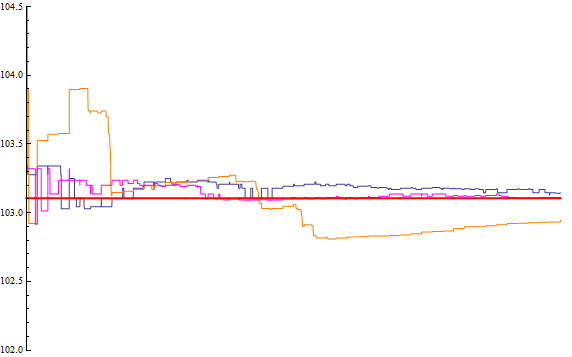}
  \caption{$1^{st}$ percentile level of the portfolio value of a homogeneous portfolio of size $500$ and up to $5.000$ simulations, MC (blue curve), QMC (orange curve) and hybrid simulation (pink curve).}
  \label{PercHom5000}
\end{figure}

The second showcase SC II: $s= 500, N=5.000$ concerns a portfolio of dimension 500 for which we carry out 5.000 simulations in the different methods. Figures 
\ref{ErwWertHom5000}, \ref{VolaHom5000} and \ref{PercHom5000} illustrate the case of a homogeneous portfolio. In all three figures we can see that the hybrid method performs very well and better than the other two methods. 

\begin{figure}[!ht]
  \centering
   \includegraphics[width=0.5\textwidth]{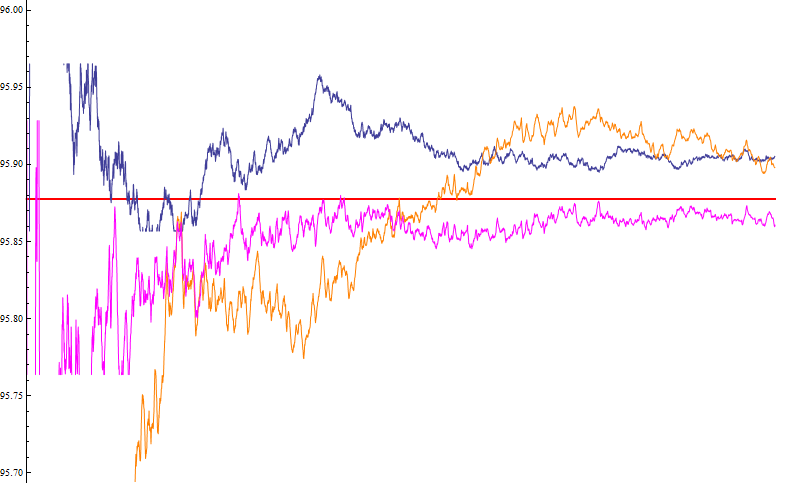}
  \caption{Expected value in the case of a size $500$ inhomogeneous portfolio and up to $5.000$ simulations, MC (blue curve), QMC (orange curve) and hybrid simulation (pink curve).}
  \label{ErwWertInHom5000}
\end{figure}

\begin{figure}[!ht]
  \centering
   \includegraphics[width=0.5\textwidth]{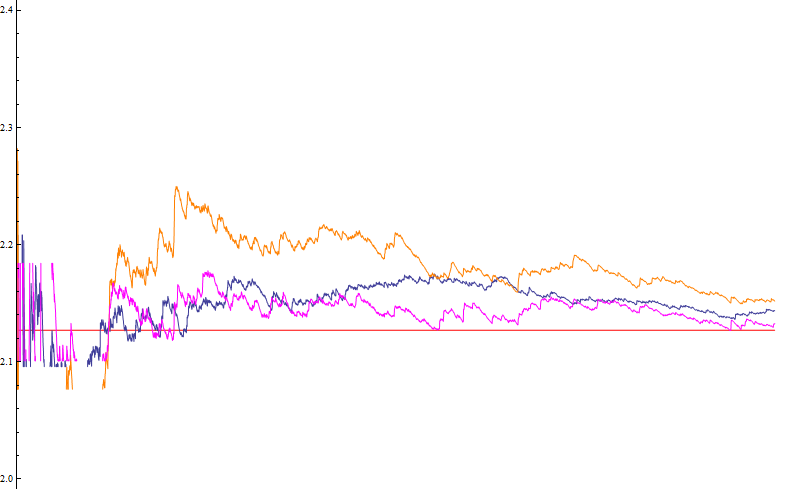}
  \caption{Standard deviation of the portfolio expected value in the case of a size $500$ inhomogeneous portfolio and up to $5.000$ simulations, MC (blue curve), QMC (orange curve) and hybrid simulation (pink curve).}
  \label{VolaInHom5000}
\end{figure}

\begin{figure}[!ht]
  \centering
   \includegraphics[width=0.5\textwidth]{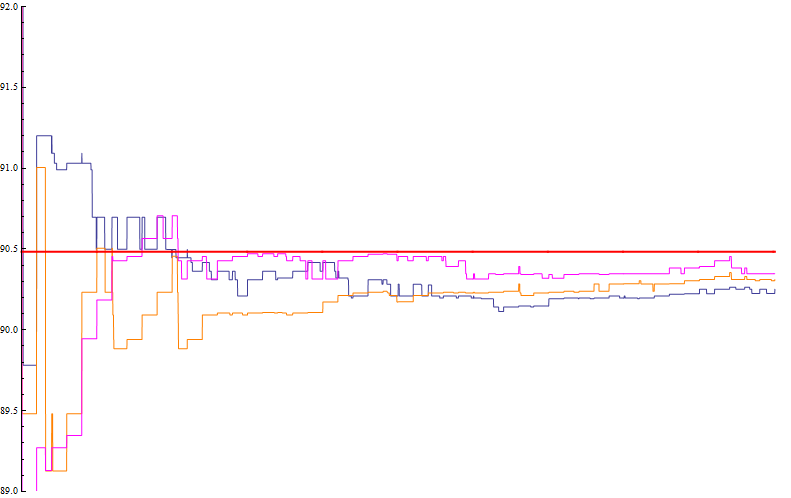}
  \caption{$1^{st}$ percentile level of the portfolio value of a inhomogeneous portfolio of size $500$ and up to $5.000$ simulations, MC (blue curve), QMC (orange curve) and hybrid simulation (pink curve).}
  \label{PercInHom5000}
\end{figure}

We compare these figures with figures \ref{ErwWertInHom5000}, \ref{VolaInHom5000}, \ref{PercInHom5000} of the example of an inhomogeneous portfolio of the same size. 
In Figure \ref{ErwWertInHom5000}  we can see again that the hybrid method performs very similarly as the MC method, both better than the QMC method.  
In Figure \ref{VolaInHom5000} and Figure \ref{PercInHom5000} we find the same better behaviour of the hybrid method in comparison with the MC and QMC methods.

\begin{figure}[!ht]
  \centering
   \includegraphics[width=0.5\textwidth]{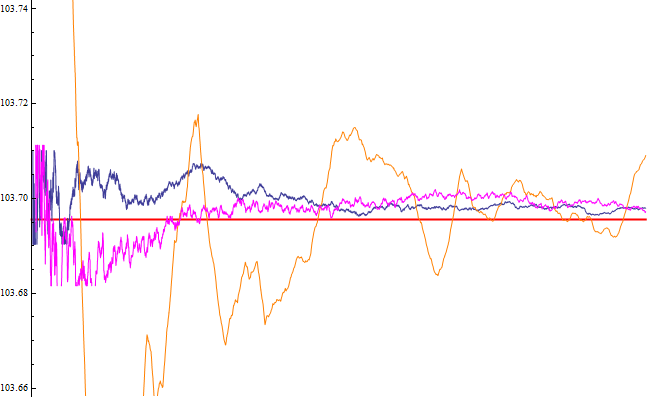}
  \caption{Expected value in the case of a size $1.000$ homogeneous portfolio and up to $10.000$ simulations, MC (blue curve), QMC (orange curve) and hybrid simulation (pink curve).}
  \label{ErwWertHom10000}
\end{figure}

\begin{figure}[!ht]
  \centering
   \includegraphics[width=0.5\textwidth]{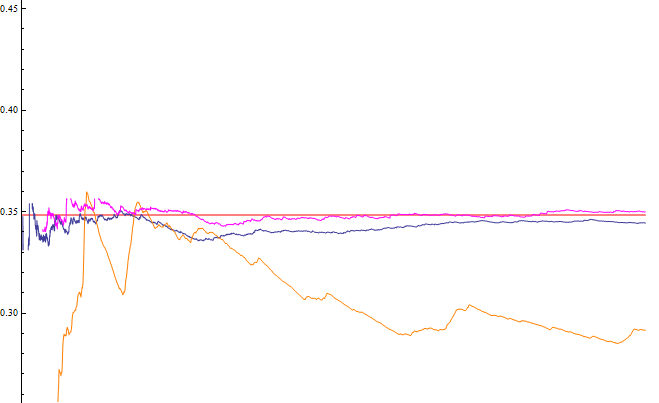}
  \caption{Standard deviation of the portfolio expected value in the case of a size $1.000$ homogeneous portfolio and up to $10.000$ simulations, MC (blue curve), QMC (orange curve) and hybrid simulation (pink curve).}
  \label{VolaHom10000}
\end{figure}

\begin{figure}[!ht]
  \centering
   \includegraphics[width=0.5\textwidth]{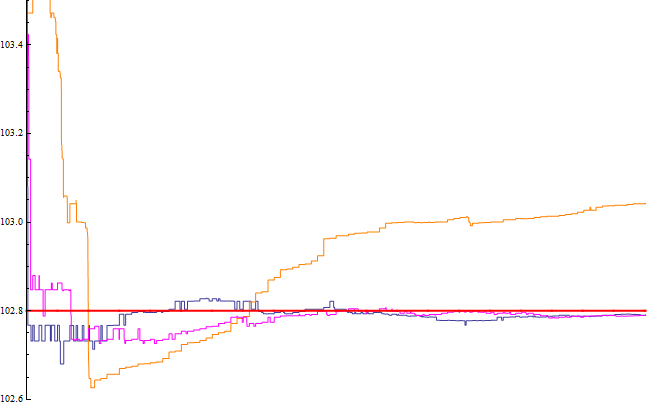}
  \caption{$1^{st}$ percentile level of the portfolio value of a homogeneous portfolio of size $1.000$ and up to $10.000$ simulations, MC (blue curve), QMC (orange curve) and hybrid simulation (pink curve).}
  \label{PercHom10000}
\end{figure}

The third and last showcase SC III: $s= 1.000, N=10.000$ concerns a portfolio of dimension 1.000 for which we carry out 10.000 simulations for the different methods. Figure \ref{ErwWertHom10000}, \ref{VolaHom10000} and \ref{PercHom10000} show the case of a homogeneous portfolio. In all three figures we can see that the hybrid method performs slightly better than the MC method. 

\begin{figure}[!ht]
  \centering
   \includegraphics[width=0.5\textwidth]{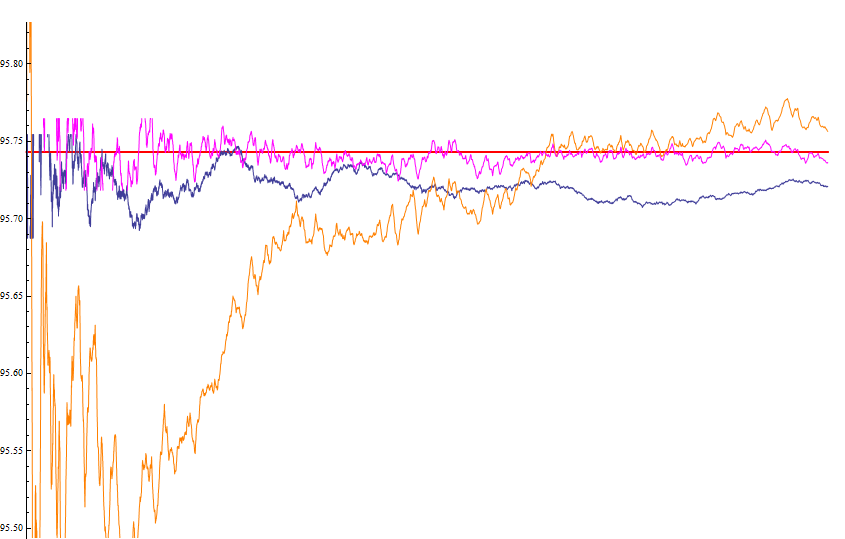}
  \caption{Expected value in the case of a size $1.000$ inhomogeneous portfolio and up to $10.000$ simulations, MC (blue curve), QMC (orange curve) and hybrid simulation (pink curve).}
  \label{ErwWertInHom10000}
\end{figure}

\begin{figure}[!ht]
  \centering
   \includegraphics[width=0.5\textwidth]{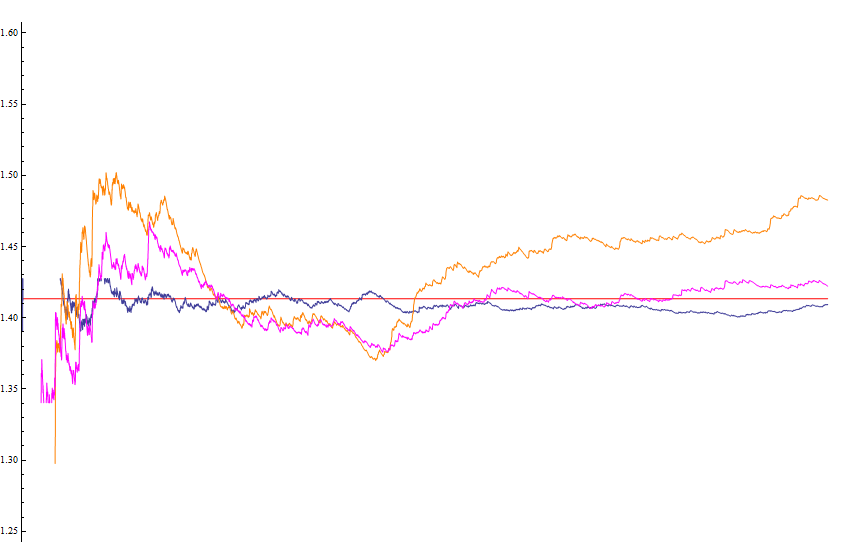}
  \caption{Standard deviation of the portfolio expected value in the case of a size $1.000$ inhomogeneous portfolio and up to $10.000$ simulations, MC (blue curve), QMC (orange curve) and hybrid simulation (pink curve).}
  \label{VolaInHom10000}
\end{figure}

\begin{figure}[!ht]
  \centering
   \includegraphics[width=0.5\textwidth]{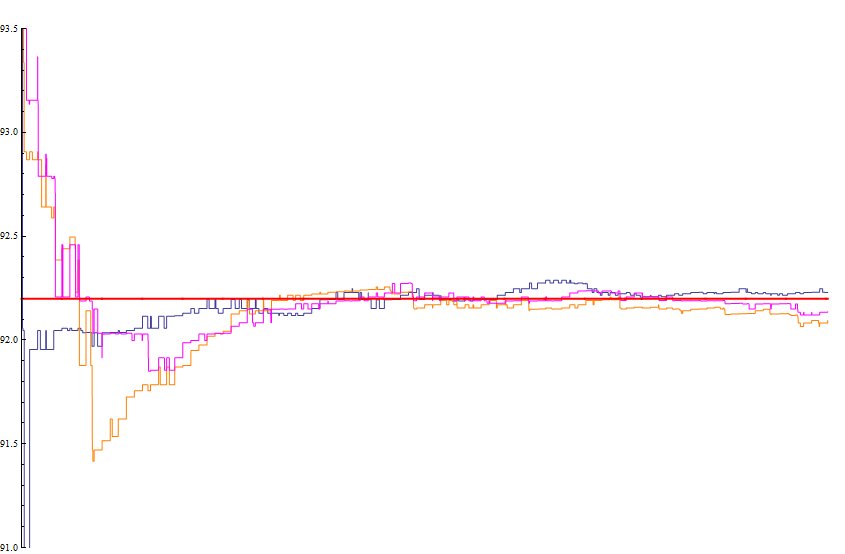}
  \caption{$1^{st}$ percentile level of the portfolio value of a inhomogeneous portfolio of size $1.000$ and up to $10.000$ simulations, MC (blue curve), QMC (orange curve) and hybrid simulation (pink curve).}
  \label{PercInHom10000}
\end{figure}

Finally we compare figures \ref{ErwWertHom10000}, \ref{VolaHom10000} and \ref{PercHom10000} with figures \ref{ErwWertInHom10000}, \ref{VolaInHom10000} and \ref{PercInHom10000} of the case of an inhomogeneous portfolio of the same size. In all $3$ figures again we observe that the hybrid method performs better than the QMC and even slightly better than the MC method. 

For us it is a bit surprising that the outperformance of the hybrid method compared with the pure MC and with QMC is not significantly better in the inhomogeneous cases than in the homogeneous cases as one would expect by our reasoning in Chapter 1.

\section{Conclusion}
 
The aim of this paper is to illustrate some examples of the quite natural use of a hybrid-Monte Carlo method for the credit risk management system CreditMetrics. Due to the structure of the problem itself, credit risk management seems to be a ideal field for the use of hybrid methods. The results obtained in this article legitimate a further and deeper investigation of these methods in applications of finance in general and of credit risk management in particular. The results of our simulations for CreditMetrics show that hybrid simulation methods perform very well in all considered cases: as expected in most cases considerably better than QMC and often better (at least not worse) than MC methods.  

However, even in the special case of CreditMetrics a lot of further work is required to make this analysis more complete. 
First of all, one should incorporate correlations between the credits in our analysis. 
Second, it would be very interesting to develop strategies to make the number of the credit products simulated with the QMC part of the hybrid sequence variable depending on the size of the portfolio, the number of the simulations, and the risk weight of the credit products in the portfolio.  
Furthermore one could investigate if other QMC sequences and other hybrid sequences can be more appropriate to be used in this context than the Niederreiter sequences.

\section*{Acknowledgment}

The authors would like to thank Isabel Pirsic for some fruitful discussions about implementations of very high dimensional Niederreiter point set.

\newpage

\section*{Figures}

\newpage

\begin{thebibliography}{5}

 
\bibitem{Spanier} Spanier, J.(1995). Quasi-Monte Carlo Methods for Particle Transport Problems. In:
\emph{Monte Carlo and Quasi-Monte Carlo Methods in Scientific Computing (H.Niederreiter and P.J.-S.Shiue, eds), Lecture Notes in Statistics}, Springer, New York, Vol. 106, p.p. 121-148.
\bibitem{Keller} Keller, A.(2013). Quasi-Monte Carlo Image Synthesis in a Nutshell. In:
\emph{Monte Carlo and Quasi-Monte Carlo Methods 2012 (J.Dick, F.Y. Kuo, G.W. Peters, and I.H. Sloan, eds)}, Springer-Verlag, Berlin, Heidelberg, Vol. 65, p.p. 2013-249.
\bibitem{Rosi} Gómez-Pérez, D., Hofer, R., Niederreiter, H.(2013). A General Discrepancy Bound for Hybrid Sequences Involving Halton Sequences. In:
\emph{Uniform Distribution Theory 8}, no.1, p.p. 31-45.
\bibitem{Oekten} Ökten, G.(1998). Applications of a Hybrid-Monte Carlo Sequence to Option Pricing. In:
\emph{Monte Carlo and Quasi-Monte Carlo Methods 1998 (H.Niederreiter and J.Spanier, eds)}, Springer-Verlag, Berlin, Heidelberg, New York, p.p. 391 -406. 
\bibitem{Nieder} Niederreiter, H.(1992). Random Number Generation and Quasi-Monte Carlo Methods. SIAM, Philadelphia.
\bibitem{Morgan} Morgan, J.P.(1997). CreditMetrics$^{\mbox{\footnotesize{TM}}}$-Technical Document. New York.
\bibitem{Nieder1} Niederreiter, H.(2009). On the discrepancy of some hybrid sequences. In:
\emph{Acta Arith. 138}, p.p. 373-398.
\bibitem{Nieder2} Niederreiter, H.(2010). A discrepancy bound for hybrid sequences involving digital explicit inversive pseudorandom numbers. In:
\emph{Unif. Distrib. Theory 5}, p.p. 53-63.
\bibitem{Nieder3} Niederreiter, H.(2010). Further discrepancy bounds and an Erd\H{o}s-Tur\'an-Koksma inequality for hybrid sequences. In:
\emph{Monatsh. Math. 161}, p.p. 193-222. 
\bibitem{Nieder4} Niederreiter, H.(2012). Improved discrepancy bounds for hybrid sequences involving Halton sequences. In:
\emph{Acta Arith. 155}, p.p. 71-84. 
\bibitem{Winterhof} Niederreiter, H., Winterhof, A.(2011). Discrepancy bounds for hybrid sequences involving digital explicit inversive pseudorandom numbers. In:
\emph{Unif. Distrib. Theory 6}., no.1, p.p. 35-56.
\bibitem{Pill} Dick, J., Pillichshammer, F.(2010). Digital Nets and Sequences. Discrepancy Theory and Quasi–Monte Carlo Integration.
Cambridge University Press. 


\end{thebibliography}
\end{document}